\documentclass[12pt,reqno]{amsart}

\usepackage{amsmath,amssymb,amscd,amsthm,bm}
\usepackage[all]{xy}
\usepackage[latin1]{inputenc}
\usepackage[T1]{fontenc}
\usepackage[hidelinks]{hyperref}
\usepackage{ifthen}
\usepackage{fullpage}
\usepackage{appendix}
\usepackage{url}
\usepackage{enumitem}
\usepackage{xcolor}
\usepackage{tikz}

\newboolean{workmode}
\setboolean{workmode}{false}

\newboolean{sections}
\setboolean{sections}{true}

\newcommand{\NN}{{\mathbb{N}}}  
\newcommand{\N}{{\mathbb{N}}}  
\newcommand{\RR}{{\mathbb{R}}}  
\newcommand{\R}{{\mathbb{R}}}  

\newcommand{\supp}{{\operatorname{supp}}}  

\newcommand{\CIN}{{C^\infty}}   

\newcommand{\Di}{{\mathfrak{Di}}} 
\newcommand{\Ch}{{\mathfrak{Ch}}} 
\newcommand{\Diffeol}{{\mathbf{Diffeol}}} 
\newcommand{\Chen}{{\mathbf{Chen}}} 

\newcommand{\ifwork}[1]{\ifthenelse{\boolean{workmode}}{#1}{}}
\newcommand{\comment}[1]{}
\newcommand{\mute}[2]{}
\newcommand{\printname}[1]{}

\ifwork{
  \renewcommand{\comment}[1]{{             \ \scriptsize{#1}\ }}
\renewcommand{\mute}[2]{{\scriptsize \ #1\ }
}
\renewcommand{\printname}[1]
    {{\makebox[0pt]{\hspace{-1.0in}\raisebox{8pt}{\tiny #1}}}}
}

\newcommand{\labell}[1] { \printname{#1} \label{#1} }

\newcommand{\ifsection}[2]{\ifthenelse{\boolean{sections}}{#1}{#2}}

\numberwithin{equation}{section}

\theoremstyle{plain}

\newtheorem{theorem}[equation]{Theorem}

\newtheorem{proposition}[equation]{Proposition}

\newtheorem{corollary}[equation]{Corollary}

\newtheorem{lemma}[equation]{Lemma}
\newtheorem*{lemma*}{Lemma}
\newtheorem*{Construction*}{Construction}

\theoremstyle{definition}

\newtheorem*{definition*}{Definition}
\newtheorem{example}[equation]{Example}
\newtheorem*{example*}{Example}
\newtheorem{examples}[equation]{Examples}

\newtheorem{remark}[equation]{Remark}
\newtheorem*{remark*}{Remark}

\definecolor{jaw}{rgb}{0,.5,0}  

\def \calC {{\mathcal C}}
\def \calF {{\mathcal F}}

\def \calD {{\mathcal D}}

\def \eps {\epsilon}

\def \ol  {\overline}

\makeatletter
\def\th@plain{%
  \thm@notefont{}
  \itshape 
}
\def\th@definition{%
  \thm@notefont{}
  \normalfont 
}
\makeatother

\title{Smooth maps on convex sets}

\date{\today}

\author{Yael Karshon}
\address{Department of Mathematics, University of Toronto,
Toronto, Ontario, Canada; and: School of Mathematical Sciences,
Tel-Aviv University, Tel-Aviv, Israel}
\email{karshon@math.toronto.edu, yaelkarshon@tauex.tau.ac.il}

\author{Jordan Watts}
\address{Department of Mathematics, Central Michigan University, Mount Pleasant, Michigan, USA 48859}
\email{jordan.watts@cmich.edu}

\begin{document}

\thanks{2020 Mathematics Subject Classification: Primary 57R55; Secondary 58-02}

\keywords{convex, diffeology, Sikorski structure,
Fr\"olicher structure, Chen structure}

\maketitle

\date{\today}

\begin{abstract}
There are several notions of a smooth map 
from a convex set to a cartesian space.  
Some of these notions coincide, but not all of them do.
We construct a real-valued function on a convex subset of the plane
that does not extend to a smooth function on any open neighbourhood
of the convex set, but that for each $k$ extends to a $C^k$ function
on an open neighbourhood of the convex set.
It follows that the diffeological and Sikorski notions of smoothness
on convex sets do not coincide.
We show that, for a convex set that is locally closed,
these notions do coincide.
With the diffeological notion of smoothness for convex sets, 
we then show that the category of diffeological spaces
is isomorphic to the category of so-called exhaustive Chen spaces.
\end{abstract}

\section{Introduction}\labell{s:intro}

There are many generalisations of smooth manifolds in the literature.
Two of them are diffeological spaces and Chen spaces.  A diffeology
on a set $X$ is given by a family of maps (called ``plots'') 
from open subsets of cartesian spaces to $X$; 
a Chen structure on a set $X$ is given by a family of maps
(called ``Chen plots'') from convex subsets of cartesian spaces to $X$.
Every Chen space has an underlying diffeology, and one can ask whether
the Chen structure can be recovered from this diffeology.  It cannot.
In fact, Stacey proves that the categories of diffeological spaces and
of Chen spaces are not equivalent \cite[Corollary 8.9]{stacey}.

The problem is that different Chen structures on the same set may have the
same underlying diffeology. 
Even the standard diffeology on the closed interval $[0,1]$ 
can come from different Chen structures
(see Example~\ref{nonstandard interval}, which is due to Stacey). 
However, we show that among the Chen structures 
with the same underlying diffeology there is a maximal one, 
characterised by a property that we call ``exhaustive'', 
and after restricting to exhaustive Chen spaces, 
the natural functor from Chen spaces to diffeological spaces
becomes an isomorphism of categories (Theorem~\ref{t:equiv}).

The details behind the discussion above, and moreover the mere definition 
of a Chen space, rely on the notion of a smooth map from a convex
set to a cartesian space.  Chen does not define this notion.
We are aware of several such notions in the literature, 
which all coincide when the convex set is closed or locally closed 
(\emph{e.g.}\ see Remark~\ref{r:locally convex} 
and Proposition~\ref{extend to boundary}).
However, these notions 
do not all coincide when the convex set is not locally closed. 
We demonstrate this by constructing 
a real-valued function on a convex subset $X$ of $\R^2$ 
that, for each $k$, 
extends to a $C^k$ function on an open neighbourhood of $X$,
but that does not extend to a $\CIN$ function on any open neighbourhood 
of $X$ (Proposition \ref{p:f}).

Thus, one needs to be careful when considering the definition
of a Chen space.  The applications that appear in Chen's papers
\cite{chen1973,chen1975,chen1977,chen1986}, however, as well as other
important works such as those of Stacey \cite{stacey} and Baez-Hoffnung
\cite{BH}, either explicitly or implicitly take the diffeological
notion of smoothness or one equivalent to it.
Moreover, this is
the notion of smoothness for which we have the isomorphism 
of categories between exhaustive Chen spaces and diffeological spaces.

Another generalization of smooth manifolds is Sikorski spaces.
A Sikorski structure on a set $X$ is given by a family of 
real-valued functions on $X$.
For a real-valued function on a convex set $X$,
the Sikorski notion of smoothness requires local extensions to the  
ambient cartesian space, whereas the diffeological notion
coincides with a more classical notion requiring
the partial derivatives on the interior to extend continuously
to the boundary (see Proposition~\ref{extend to boundary},
which relies on a result of Kriegl).

The best possible scenario is when there are both a diffeology
and a Sikorski structure that determine each other.
We refer to such structures as \emph{reflexive}.
Such structures form a category that is isomorphic to the category of Fr\"olicher spaces; see \cite{BIZKW} and \cite[Chapter 2]{watts}.  For more details on Fr\"olicher spaces, see \cite{frolicher}.

Arbitrary convex sets are not necessarily reflexive;
this is our Corollary \ref{cor:not reflexive},
which follows from our Proposition \ref{p:f}.
In contrast, closed convex sets \emph{are} reflexive;
this is our Theorem~\ref{t:convex}, whose proof relies on a result of Kriegl.
Moreover, 
in Section \ref{sec:closed} we show that reflexivity is a local property
(with a careful interpretation of ``local''; 
see Propositions~\ref{p:local reflexivity} and~\ref{p:open subset}).
As a result, we conclude that any Sikorski space that is locally 
Sikorski diffeomorphic to closed convex subsets of cartesian spaces
is reflexive (Theorem~\ref{t:locally convex}).

And so, this paper has three purposes.  The first is to demonstrate
that various notions of smoothness 
(specifically, the diffeological and Sikorski notions)
do not coincide on convex sets that are not locally closed.
We achieve this in Section~\ref{s:example}.
The second is to show that Sikorski spaces 
that are locally diffeomorphic to closed convex sets,
such as manifolds-with-corners, are reflexive.
We achieve this in Sections~\ref{sec:closed} and~\ref{s:locality},
after establishing ``locality of reflexivity'',
which is interesting in its own right.
The third is to show the isomorphism of categories 
between diffeological spaces and exhaustive Chen spaces. 
We achieve this in Section \ref{sec:chen}.

\subsection*{Acknowledgements}
This paper is motivated by correspondence between the authors and Patrick
Iglesias-Zemmour many years ago on the question of whether the category
of Chen spaces is equivalent to the category of diffeological spaces.

Yael Karshon acknowledges the support
of the Natural Sciences and Engineering Research Council of Canada
(NSERC).

\section{Preliminaries}\labell{s:prelim}

We recall the definitions of a Sikorski structure
and of a diffeology.  
Good references are, respectively, \'Sniatycki's book \cite{sniatycki}
and Iglesias-Zemmour's book \cite{iglesias}.  Note that in the literature, including \'Sniatcyki's book and our paper \cite{BIZKW}, Sikorski structures are often called ``differential structures''.

A \textbf{cartesian space} is $\R^d$ for some $d \in \N$, 
equipped with its standard smooth structure.
A \textbf{convex set} is a convex subset of a cartesian space.

Let $X$ be a set.
A \textbf{Sikorski structure} on $X$ is a non-empty 
collection $\calF_X$ of real-valued functions on $X$
that satisfies the following two axioms.
We equip $X$ with the \textbf{initial topology} induced by~$\calF_X$,
which is the weakest topology for which all the functions in $\calF_X$
are continuous. 
\begin{itemize}
\item[(S1)]
If $f_1,\ldots,f_n \in \calF_X$ and $g \in \CIN(\R^n)$,
then $g \circ (f_1,\ldots,f_n) \in \calF_X$.
\item[(S2)]
Given $f \colon X \to \R$, if for each point of $X$
there exist a neighbourhood $V$ and a function $g \in \calF_X$
such that $g|_V = f|_V$, then $f \in \calF_X$.
\end{itemize}

A \textbf{diffeology} on $X$ is a collection $\calD_X$
of maps from open subsets of cartesian spaces to $X$,
called plots,
that satisfies the following three axioms.
\begin{itemize}
\item[(D1)]
Constant maps to $X$ are plots.
\item[(D2)]
The precomposition of a plot with a $\CIN$ map
between open subsets of cartesian spaces is a plot. 
\item[(D3)]
Given $p \colon U \to X$, if each point of $U$
has a neighbourhood $V$ such that $p|_V$ is a plot,
then $p$ is a plot.
\end{itemize}

Any manifold $X$ has a natural Sikorski structure and a natural diffeology,
consisting of those functions $X \to \R$ and those maps $U \to X$
that are $\CIN$ in the usual sense.

Given Sikorski spaces $(X,\calF_X)$ and $(Y,\calF_Y)$,
a map $X \to Y$ is \textbf{Sikorski smooth}
if its composition with every function in $\calF_Y$ is in $\calF_X$.
Given diffeological spaces $(X,\calD_X)$ and $(Y,\calD_Y)$,
a map $X \to Y$ is \textbf{diffeologically smooth}
if its precomposition with every plot of $X$ is a plot of $Y$.  
Both Sikorski spaces with Sikorski smooth maps between them, 
and diffeological spaces with diffeologically smooth maps between them, 
form categories, denoted $\mathbf{Sik}$ and $\Diffeol$, resp.  
Both categories contain the category of smooth manifolds as full 
subcategories.

When considering diffeologies and Sikorski structures
on the same set $X$, it is convenient to combine them into one structure.
This approach was introduced in the introduction of 
Watts' thesis \cite{watts};
its topological aspects are studied in Virgin's paper \cite{virgin}.  
Thus, we consider triples $(\mathcal{D}_X,X,\mathcal{F}_X)$,
where $\calD_X$ is a diffeology on $X$
and $\calF_X$ is a Sikorski structure on $X$,
that are \textbf{compatible} in the following sense:
for any $p \colon U \to X$ in $\calD_X$ and $f \colon X \to \R$ 
in $\calF_X$, the composition $f\circ p \colon U \to \R$ is smooth. 

Given a Sikorski space $(X,\calF_X)$,
the set of Sikorski smooth maps
from open subsets of cartesian spaces to $X$
is a diffeology, which we denote $\Pi\calF_X$.
Given a diffeological space $(X,\calD_X)$,
the set of diffeologically smooth real-valued functions on $X$
is a Sikorski structure, which we denote $\Phi\calD_X$.
A Sikorski structure $\calF_X$ is \textbf{reflexive}
if $\Phi\Pi\calF_X = \calF_X$.
A diffeology $\calD_X$ is \textbf{reflexive}
if $\Pi\Phi\calD_X = \calD_X$.
For any Sikorski structure $\calF_X$,
the diffeology $\Pi\calF_X$ is reflexive.
For any diffeology $\calD_X$,
the Sikorski structure $\Phi\calD_X$ is reflexive.
A triple is $(\calD_X, X, \calF_X)$ is \textbf{reflexive} if
$\Phi\mathcal{D}_X=\mathcal{F}_X$ and $\Pi\mathcal{F}_X=\mathcal{D}_X$.

For any Sikorski space $(X,\calF_X)$ and subset $V$ of $X$,
the \textbf{subspace Sikorski structure} 
is the set of those functions $f \colon V \to \R$ such that 
for each point of~$V$ there is a neighbourhood $W$ in $X$
and a function $g$ in $\calF_X$
such that $g|_{V \cap W} = f|_{V \cap W}$.
For any diffeological space $(X,\calD_X)$ and subset $V$ of $X$,
the \textbf{subset diffeology} 
is the set of those maps $p \colon U \to V$
from open subsets $U$ of cartesian spaces to $V$
whose composition with the inclusion map $V \to X$ is in $\calD_X$.

For any subset $X$ of any cartesian space,
except where we say otherwise,
we equip $X$ with its subspace Sikorski structure $\calF_X$
and with its subset diffeology $\calD_X$
that are induced from the ambient cartesian space.
When $X$ is open in its ambient cartesian space, 
these structures coincide with its Sikorski and diffeological structures 
as a manifold.

The \textbf{D-topology} on a diffeological space
is the strongest topology for which all the plots are continuous.

{
\begin{lemma} \labell{subsets}
For any subset $X$ of any cartesian space $\R^d$,
\begin{enumerate}[itemsep=6pt]
\item $\calD_X = \Pi\calF_X$;
\item The initial topology on $X$ coincides with its subspace topology;
\item
If $X$ is convex, then the D-topology on $X$
coincides with its subspace topology.
\end{enumerate}
\end{lemma}

\begin{proof}
Item (1) is straightforward. Item (2) is \cite[Lemma 2.28]{watts}.  
Item (3) is \cite[Theorem 3.17]{CSW}.
\end{proof}
}

\section{An example}\labell{s:example}

Let $X$ be the union of the open upper half plane
and the non-negative real axis in $\RR^2$:

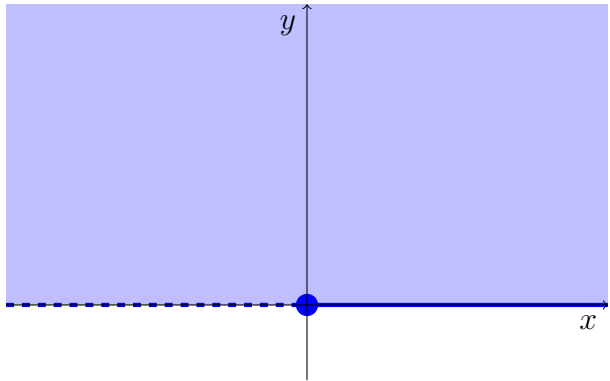
\begin{center}
\begin{tikzpicture}
		\fill[blue!25] (-4,0) rectangle (4,4);
		\draw[blue, ultra thick, dashed]  (-4,0) -- (0,0);
		\draw[blue, ultra thick] (0,0) -- (4,0);
		\filldraw[blue] (0,0) circle (4pt);
		\draw[->] (-4,0) -- (4,0) node[below left] {$x$};
		\draw[->] (0,-1) -- (0,4) node[below left] {$y$};
\end{tikzpicture}
\end{center}

$$ X :=  \{ (x,y) \in \R^2 \ | \ y > 0 \} \ \cup \ 
         \{ (x,0) \ | \ x \geq 0 \} .$$

\begin{proposition} \labell{p:f}
There exists a function 
$$ f \colon X \to \R $$
that does not extend to a $\CIN$ function 
on any open neighbourhood of $X$ in $\R^2$,
but that, for each $k \in \N$, extends to a $C^k$ function 
on some open neighbourhood of $X$ in $\R^2$.  
\end{proposition}

\begin{proof}
For each $m\in \NN$, let $\varphi_m \colon \RR \to \R$ be a $C^m$ function that vanishes on $(-\infty,0]$, is smooth on $(0,\infty)$, whose derivatives of all orders less than or equal to $m$ are bounded, and whose $m+1^\text{st}$ derivative is not bounded on any punctured neighbourhood of $0$.  For example, we can take such a $\varphi_m$ that vanishes outside of $(0,2)$ and is equal to $x^{m+1/2}$ on $(0,1)$.

Let $b\colon\RR\to\RR$ be a smooth function that is equal to $0$ on $(-\infty,-1]$, is equal to $-1$ on $[0,\infty)$, and is strictly decreasing on $[-1,0]$.  For each $0<\eps\leq 1$, let $b_\eps(x):=\eps b(x-1+\eps)$.  Note that $\eps\mapsto b_\eps(x)$ is decreasing for all $x$ and strictly decreasing for all $x\geq 0$.

For each $m\in\NN$, define $h_m\colon\RR^2\to\RR$ by $$h_m(x,y):=y-b_{1/m}(x).$$  Then $h_m$ is smooth, and every iterated derivative of $h_m$ of every positive order is bounded on $\RR^2$.  Let $$U_m:=\{(x,y)\mid h_m(x,y)>0\}.$$  Then $U_1\supset U_2\supset\dots$, and $X$ is the intersection of the open sets $U_m$.

Because the derivatives of order less than or equal to $m$ of $\varphi_m$ and of $h_m$
are bounded, so are the derivatives of order less than or equal to $m$
of $\varphi_m \circ h_m$. Let $c_m>0$ be
a bound on the absolute values of these derivatives.  
Define
$$f := \sum_{m=1}^\infty \frac{1}{c_m 2^m} \varphi_m \circ h_m\Big|_X.$$  

For each $k\in\NN$, the series $$g_{k+1}:=\sum_{m=k+1}^\infty\frac{1}{c_m2^m}\varphi_m\circ h_m$$ converges to a $C^{k+1}$ function $\RR^2\to\RR$.  For $1\leq m\leq k$, the functions $\varphi_m\circ h_m$ are smooth on the set $U_k$.  Thus the series $$\sum_{m=1}^\infty\frac{1}{c_m2^m}\varphi_m\circ h_m\Big|_{U_k}$$ converges to a $C^{k+1}$ function on $U_k$ whose restriction to $X$ is $f$.

	We claim that $f$ does not extend to a smooth function on any neighbourhood of the origin.  Indeed, assuming otherwise, fix $k\in \NN$ such that $(-1/k,0)$ is in such a neighbourhood.  For $1\leq m < k$, we have $h_m(-1/k,0)>0$, and so the functions $\varphi_m\circ h_m$ are smooth on some neighbourhood of $(-1/k,0)$.  Since $g_{k+1}$ is $C^{k+1}$ on $\RR^2$, the series $$\sum_{\substack{m=1 \\ m \neq k}}^\infty\frac{1}{c_m2^m}\varphi_m\circ h_m$$ converges to a $C^{k+1}$ function on some neighbourhood of $(-1/k,0)$.  Subtracting it from $f$, we obtain that $\varphi_k\circ h_k|_X$ extends to a $C^{k+1}$ function on some neighbourhood of $(-1/k,0)$.  This implies that the $k+1^\text{st}$ derivatives of $\varphi_k\circ h_k$ are bounded on the intersection of the interior of $X$ with some neighbourhood of $(-1/k,0)$.  But $\varphi_k(h_k(-1/k,y))=\varphi_k(y)$, and so the $k+1^\text{st}$ derivative of $\varphi_k\circ h_k$ with respect to $y$ is not bounded on the intersection of $\{-1/k\}\times(0,\infty)$ with any neighbourhood of $(-1/k,0)$.
\end{proof}

\begin{corollary} \labell{cor:not reflexive}
The Sikorski structure on $X$ is not reflexive.
\end{corollary}

\begin{proof}
Let $f \colon X \to \R$ be a function as in Proposition \ref{p:f}.
Then $f$ is diffeologically smooth.
Indeed, let $p \colon U \to X$ be a plot, \emph{i.e.}, a map
from an open subset $U$ of some cartesian space $\R^\ell$ to~$X$
that is smooth as a map to $\R^2$.
For each $k \in \N$, because $f$ extends to a $C^k$ function 
on some neighbourhood of $X$ in $\R^2$, 
the composition $f \circ p \colon U \to \R$ is $C^k$.
Because this holds for all $k \in \N$,
the composition $f \circ p \colon U \to \R$ is smooth.
Because the plot $p$ is arbitrary, 
$f \colon X \to \R$ is diffeologically smooth.

However, $f$ is not Sikorski smooth, 
because it does not extend to a $C^\infty$ function 
on any open neighbourhood of $X$ in $\R^2$.
Because the subset diffeology on $X$ in induced from the 
subspace Sikorski structure on $X$ (Lemma \ref{subsets}),
the Sikorski structure on $X$ is not reflexive.
\end{proof}

\section{Closed convex sets}
\labell{sec:closed}

\nopagebreak

The phenomenon exemplified in Corollary~\ref{cor:not reflexive}
does not occur with convex subsets of $\R^n$ that are closed:

\begin{theorem}[Closed Convex Sets are Reflexive]\labell{t:convex}
Let $K\subseteq\RR^n$ be closed and convex.
Then the Sikorski structure on $K$ is reflexive.
\end{theorem}

In preparation for the proof of Theorem~\ref{t:convex},
we use results of Boman and Kriegl to show the following result.

\begin{proposition} \labell{extend to boundary}
Let $K$ be a convex subset of $\R^n$ with non-empty interior $\mathring{K}$.
A map $f\colon K\to \RR$ is diffeologically smooth if and only if 
$f|_{\mathring{K}} \colon \mathring{K} \to \R$ is smooth
and all of its partial derivatives of all orders extend continuously to $K$.
\end{proposition}

\begin{proof}
We will use a result of Kriegl.
Kriegl allows convex subsets in infinite dimensions.
He works with the $c^\infty$-topology,
and with a notion of smoothness to which we will refer here
as $c^\infty$-smoothness.
These notions rely on the notion of a smooth curve in $K$.
When $K$ is a subset of $\RR^n$,
a smooth curve in $K$ is a map $c \colon \RR \to K$
that is smooth (infinitely-differentiable) in the usual sense 
as a map to $\R^n$.

The $c^\infty$-topology on $K$ is defined to be the strongest topology on $K$ 
for which all smooth curves $c \colon \R \to K$ are continuous.  
In our finite dimensional setup, the $c^\infty$-topology coincides
with the subspace topology.
Indeed, the $c^\infty$-topology coincides with the D-topology
by \cite[Theorem 3.7]{CSW},
and the D-topology coincides with the standard topology
by \cite[Theorem 3.17]{CSW}.
Thus, we may refer without ambiguity to the interior $\mathring{K}$,
and to a function on $\mathring{K}$ having a continuous extension to $K$.

Suppose $f$ is diffeologically smooth. In particular, for all smooth curves
$c\colon\RR\to K$, the composition $f\circ c$ is smooth. By definition,
this means that $f$ is $c^\infty$-smooth.  Conversely, suppose $f$ is
$c^\infty$-smooth.  
Let $p\colon U\to K$ be a plot of $K$, that is, a smooth map in the
usual sense.  For any smooth curve $c\colon\RR\to U$, the composition
$p\circ c$ is a smooth curve in $K$.  Since $f$ is $c^\infty$-smooth,
the composition $f\circ p\circ c$ is smooth.  By a result of Boman
\cite[Theorem 1]{boman}, since $c$ is arbitrary, $f\circ p$ is smooth.
Since $p$ is arbitrary, $f$ is diffeologically smooth.

By a result of Kriegl \cite[Theorem 1.5]{kriegl},
$c^\infty$-smoothness of $f$ is equivalent to the restriction of $f$ to the interior $\mathring{K}$ of $K$ being smooth
and it and all of its derivatives 
extending to functions on $K$ continuously 
(with respect to the $c^\infty$-topology on $K$,
hence with respect to the standard topology on $K$).

We pause to explain Kriegl's term \emph{derivatives}.
Kriegl works with the setup of Fr\"olicher and Kriegl \cite{FK}.
Unravelling \cite[Definitions~4.1.9, 4.3.9, and~4.3.26]{FK}
for the special case of a smooth real-valued function 
on an open subset of $\R^n$,
the derivatives of $f|_{\mathring{K}} \colon \mathring{K} \to \R$ 
are the same as its iterated differentials in the usual sense.
These having continuous extensions to $K$
is equivalent to all the partial derivatives of $f$ of all orders
(in the usual sense) having continuous extensions to $K$.
\end{proof}

We are now ready to prove Theorem~\ref{t:convex}.

\begin{proof}[Proof of Theorem~\ref{t:convex}]
Recall (Lemma~\ref{subsets}) that the diffeology on $K$
is induced from the Sikorski structure on $K$.
So we need to show that, for every function $f \colon K \to \R$,
if $f$ is diffeologically smooth, then $f$ is Sikorski smooth.
By replacing $\RR^n$ with the affine span of~$K$, 
we may assume that the interior of $K$ is non-empty.
Fix a function $f \colon K \to \R$ that is diffeologically smooth.
By Proposition~\ref{extend to boundary},
the restriction of $f$ to the interior $\mathring{K}$ of $K$ is smooth
and all its partial derivatives of all orders
extend to continuous functions on~$K$.
Because $K$ is closed and convex, it satisfies the hypothesis 
of the Whitney Extension Theorem 
(see, for instance, 
\cite[Theorem 2.6 and Proposition 2.16]{bierstone:diff fcts}),
and hence $f$ extends to a smooth function $g\colon\RR^n\to\RR$.
So $f$ is Sikorski smooth.
\end{proof}

\begin{example}\labell{x:square pyramid}
Let $P$ be the solid square pyramid.  By Theorem~\ref{t:convex}, 
$P$ obtains a reflexive Sikorski structure and a reflexive
diffeological structure from $\RR^3$.  One may interpret this to mean
that $P$ has an unambiguous smooth structure,
although it is not a manifold-with-corners 
({\it cf.} Example~\ref{x:mfld-w-corners}).
\end{example}

\section{Locality of reflexivity}
\labell{s:locality}

\nopagebreak

In this section, we show that compatibility and reflexivity of a triple
$(\calD, X, \calF)$ are local properties; see Propositions~\ref{p:local reflexivity} and~\ref{p:open subset}.
We apply this to subcartesian spaces
that are locally compact and locally convex; see Theorem~\ref{t:locally convex}.
In particular, manifolds-with-corners are reflexive;
see Example~\ref{x:mfld-w-corners}. 

We begin with a lemma.

\begin{lemma} \labell{locality lemma}
Let $(X,\calF)$ be a Sikorski space,
let $(V,\calF_V)$ be a Sikorski subspace, and let $p \in \Pi\calF$
with image in $V$.  Then $p \in \Pi\calF_V$.
\end{lemma}

\begin{proof}
As usual, we equip $X$ with its initial topology that is induced by $\calF$.
Let $f_V \in \calF_V$.
For each $y \in V$ there exist an open neighbourhood $W$ of $y$
and a function $g \in \calF$ such that $f_V|_{V \cap W} = g|_{V \cap W}$. 
Because $p \in \Pi\calF$,
the preimage $p^{-1}(W)$ is open and the composition $g \circ p$ is smooth.
So $f_V \circ p|_{p^{-1}(W)} = g \circ p|_{p^{-1}(W)}$ is smooth.
Because the sets $p^{-1}(W)$ form an open cover of the domain of $p$,
the composition $f_V \circ p$ is smooth.
It follows that $p \in \Pi\calF_V$.
\end{proof}

For a diffeology $\calD$ and Sikorski structure 
$\calF$ on a set $X$ and for a subset $V$ of $X$,
we denote by $\calD_V$ and $\calF_V$ the subset diffeology
and the subspace Sikorski structure
that are induced by $X$.  

We now prove locality of ``compatible'' and ``reflexive''.

\begin{proposition}[Compatibility and Reflexivity are Local]
\labell{p:local reflexivity}
Let $\mathcal{D}$ and $\mathcal{F}$ be a diffeology and a Sikorski
structure on a set $X$. 
Denote by $\tau_\calD$ the D-topology on $X$ induced by $\calD$
and by $\tau_\calF$ the initial topology on~$X$ induced by $\calF$.
\begin{enumerate}[itemsep=6pt]
\item Suppose for every point in $X$ there is a $\tau_\calD$-open neighbourhood $V$ 
such that $(\mathcal{D}_V,V,\mathcal{F}_V)$ is a compatible triple.  
Then $(\mathcal{D},X,\mathcal{F})$ is a compatible triple.
In particular, the topology $\tau_\calD$ contains the topology $\tau_\calF$. 
\item Suppose that the topology $\tau_\calD$ contains 
the topology $\tau_\calF$
and that for every point in $X$ there is a $\tau_\calF$-open neighbourhood $V$
such that $(\mathcal{D}_V,V,\mathcal{F}_V)$ is a reflexive triple.
Then $(\mathcal{D},X,\mathcal{F})$ is a reflexive triple.
\end{enumerate}
\end{proposition}

\begin{proof}
We prove Claim (1). Fix $p\in\calD$, $f\in\calF$, and $x \in X$.
By hypothesis, there is a $\tau_D$-open neighbourhood $V$ of $x$ in $X$
such that $(\mathcal{D}_V,V,\mathcal{F}_V)$ is a compatible triple.  
The restriction $p_V := p|_{p^{-1}(V)}$ takes values in $V$,
and its domain is open in the domain of $p$, so $p_V \in \calD_V$.
Since $f_V := f|_V$ is in $\mathcal{F}_V$ and by compatibility, the composition
$f\circ p|_{p^{-1}(V)}=f_V\circ p_V$ is smooth.  
Since smoothness of $f\circ p$ is a local condition
and the sets $p^{-1}(V)$ form an open cover of the domain of $p$,
this shows that $f\circ p$ is smooth.  Claim (1) follows.

We prove Claim (2).  By Claim (1), it is enough
to show that $\Pi\mathcal{F}\subseteq\mathcal{D}$ and
$\Phi\mathcal{D}\subseteq\mathcal{F}$.  
For the first inclusion, fix $p\in\Pi\mathcal{F}$.
By the locality axiom of a diffeology, it suffices to show
that for each $x \in X$ there is a $\tau_\calD$-open neighbourhood $V$
such that $p_V := p|_{p^{-1}(V)} \in \calD_V$.
Let $x \in X$. 
By hypothesis, there is a $\tau_\calF$-open neighbourhood $V$ of $x$ 
such that $(\mathcal{D}_V,V,\mathcal{F}_V)$ is a reflexive triple. 
In particular, $\calD_V = \Pi\calF_V$.  
Since $\tau_\calF$ is contained in $\tau_\calD$,
it suffices to show that $p_V \in \Pi\calF_V$.
This, in turn, follows from Lemma~\ref{locality lemma}, since $p\in\Pi\calF$.

For the second inclusion $\Phi\mathcal{D}\subseteq\mathcal{F}$, fix
$f\in\Phi\mathcal{D}$.  By the locality axiom of a Sikorski space, 
it suffices to show that for each $x\in X$ there is a $\tau_\calF$-open 
neighbourhood $W$ of $x$ 
and a function $g\in\mathcal{F}$ such that $f|_W=g|_W$.  Fix $x\in X$.  
By hypothesis, there is a $\tau_\calF$-open neighbourhood $V$ of $x$
such that $(\mathcal{D}_V,V,\mathcal{F}_V)$ is a reflexive triple. 
In particular, $\calF_V=\Phi\calD_V$.  Define $f_V:=f|_V$.
For any $p\in\mathcal{D}_V$,
because $p$ is in $\calD$ and $f$ is in $\Phi\calD$,
the composition $f_V\circ p = f \circ p$ is smooth.  
So $f_V\in\Phi\mathcal{D}_V=\calF_V$.
By definition of $\calF_V$, there exist a $\tau_\calF$-open
neighbourhood $W\subseteq V$ of $x$ and a function $g\in\mathcal{F}$ such
that $g|_W = f_V|_W=f|_W$, proving what is required.
\end{proof}

In Part (2) of Proposition~\ref{p:local reflexivity}, 
the interaction between the two topologies 
is not superfluous.  Indeed, let
$X$ be $\RR$ and $\mathcal{D}$ be the standard diffeology on it.  Take
$\mathcal{F}$ to be \emph{all} real-valued functions.  Then the initial
topology induced by $\mathcal{F}$ is discrete, and taking singleton
open sets, we satisfy the hypotheses of both claims in the proposition.
However, compatibility (let alone reflexivity) fails:
no non-constant smooth curve in $\mathcal{D}$ is in $\Pi\mathcal{F}$.

We pause to recall that the initial topology admits bump functions.

\begin{lemma} \labell{bump}
Let $(X,\calF)$ be a Sikorski space,
let $\tau_\calF$ be the initial topology on $X$ induced by $\calF$,
let $x \in X$, and let $V$ be a $\tau_\calF$-open neighbourhood 
of $x$ in $X$.  Then there exists $\rho \in \calF$
whose $\tau_\calF$-support is contained in $V$
and that is equal to~$1$ on a $\tau_\calF$-neighbourhood of $x$ in $X$.
\end{lemma}

\begin{proof}
Because $V$ is a $\tau_\calF$-open neighbourhood of $x$,
there exist $h_1,\ldots,h_N\in\calF$ such that 
$x\in \cap_{i=1}^N h_j^{-1}((0,1))\subseteq V$.
Let $h := (h_1,\ldots,h_N)$,
and let $b\colon\RR^N\to[0,1]$ be a smooth function 
that is identically~$1$ on an open neighbourhood of $h(x)$ 
and whose support is contained in $(0,1)^N$.  
Note that $b \circ h \in \calF$,
that the $\tau_\calF$-support of $b \circ h$ is contained in $V$,
and that $(b \circ h)^{-1}(\{1\})$ 
is a $\tau_\calF$-neighbourhood of $x$ in $X$.
Define $\rho \colon X \to \RR$ 
to be equal to $b \circ h$ on $V$ and $0$ elsewhere.
Then the $\tau_\calF$-support $K$ of $\rho$ coincides with that of $b \circ h$,
and $\rho^{-1}(\{1\}) = (b \circ h)^{-1}(\{1\})$.
Finally, because $\rho$ coincides with an element of $\calF$
on each of the two sets $V$ and $X \smallsetminus K$
and these two sets form an open cover of $X$,
we conclude that $\rho \in \calF$.
\end{proof}

We now prove a converse statement to Proposition~\ref{p:local reflexivity}: 
global compatibility implies local compatibility,
and global reflexivity implies local reflexivity.

\begin{proposition}\labell{p:open subset}
Let $\calD$ and $\calF$ be a diffeology and a Sikorski structure
on a set~$X$. 
If the triple $(\calD, X, \calF)$ is compatible 
and $V \subset X$ is any subset,
then the triple $(\calD_V, V, \calF_V)$ is compatible.
If the triple $(\calD, X, \calF)$ is reflexive 
and $V \subset X$ is an open subset in the initial topology,
then the triple $(\calD_V, V, \calF_V)$ is reflexive.
\end{proposition}

\begin{proof}
Suppose that $(\calD,X,\calF)$ is compatible, and let $V \subset X$.
Fix $p\in\mathcal{D}_V$, $f\in\mathcal{F}_V$, and $x\in V$.
There exist a $\tau_\calF$-open neighbourhood $W$ of $x$ in $X$
and a function $g\in\mathcal{F}$ such that $f|_{V \cap W}=g|_{V \cap W}$.
Compatibility implies that $p^{-1}(W)$ is open in the domain of $p$
and that $g \circ p$ is smooth.
So $f\circ p|_{p^{-1}(W)}=g\circ p|_{p^{-1}(W)}$ is smooth.
Because such sets $p^{-1}(W)$ are an open cover of $V$,
the composition $f\circ p$ is smooth.  
Compatibility of $(\calD_V,V,\calF_V)$ follows.

Suppose now that $(\calD,X,\calF)$ is reflexive and $V \subset X$
is $\tau_\calF$-open.
By the previous claim, it suffices to show that $\Pi\calF_V \subset \calD_V$
and $\Phi\calD_V\subset\calF_V$. For the first inclusion,
fix $p\in\Pi\mathcal{F}_V$.
For any $f\in\mathcal{F}$, the restriction $f|_V$ is in $\mathcal{F}_V$, 
so $f\circ p=f|_V\circ p$ is smooth.  Thus, $p\in\Pi\mathcal{F}=\mathcal{D}$
with image in $V$, and so $p\in\mathcal{D}_V$.

For the second inclusion $\Phi\mathcal{D}_V\subseteq\mathcal{F}_V$, 
fix $f\in\Phi\mathcal{D}_V$.  
It suffices to show that for any $x\in V$
there exist a $\tau_\calF$-open neighbourhood $W$ of $x$ in $V$
and a function $g\in\mathcal{F}$ such that $f|_{W}=g|_{W}$.  Fix $x\in V$.  
By Lemma~\ref{bump}, there exists a function 
$\rho \in \calF$ whose $\tau_\calF$-support $\supp(\rho)$ is contained in $V$
and that is equal to~$1$ on a $\tau_\calF$-neighbourhood of $x$ in $X$.
Define $g\colon X\to\RR$ 
to be equal to the product $\rho \cdot f$ on $V$ and $0$ elsewhere.  
Fix $p\in\mathcal{D}$.  Write 
$$ (g \circ p)|_{p^{-1}(V)}
 \ = \ \rho \circ p|_{p^{-1}(V)} \ \cdot \ f \circ p|_{p^{-1}(V)}.$$
By compatibility, $p^{-1}(V)$ is open in the domain of $p$,
the first factor on the right hand side is smooth,
and since $f \in \Phi\calD_V$, 
the second factor on the right hand side is smooth.
So the composition $g \circ p$ is smooth in the open subset $p^{-1}(V)$
of the domain of $p$.
Since this composition vanishes outside $p^{-1}(\supp(\rho))$,
which is a closed subset of the domain of $p$ and is contained
in $p^{-1}(V)$, this composition is smooth.
We conclude that $g \in \Phi\calD = \calF$.
Since $f$ and $g$ coincide on the $\tau_\calF$-interior of $\rho^{-1}(1)$, 
which contains $x$, this proves what is required.
\end{proof}

We now prove that Sikorski spaces
that are locally compact and locally convex are reflexive.
We start with a technical lemma:

\begin{lemma}\labell{l:locally convex}
Let $C'$ be a subset of $\R^n$.
Suppose that $C'$ is locally compact and that $C' = C \cap U$
for some convex subset $C$ of $\R^n$ and some open subset $U$ of $\R^n$.
Then the Sikorski structure of $C'$ is reflexive.
\end{lemma}

\begin{proof}
Since $C'$ is locally compact and is contained in the open set $U$,
there is a collection $\{ \ol{B}_i \}$ of closed balls in $\R^n$
that are contained in $U$, 
whose interiors $\mathring{B}_i$ cover $C'$,
and such that each $C' \cap \ol{B}_i$ is closed in $\R^n$.
By Theorem~\ref{t:convex}, the Sikorski structure 
of each $C' \cap \ol{B}_i$ is reflexive.
By Proposition~\ref{p:open subset}, the Sikorski structure 
of each $C' \cap \mathring{B}_i$ is reflexive.
By Proposition~\ref{p:local reflexivity}, 
the Sikorski structure of $C'$ is reflexive.
\end{proof}

We continue to general Sikorski spaces
that are locally compact and locally convex:

\begin{theorem} \labell{t:locally convex}
Let $X$ be a locally compact Sikorski space 
satisfying the following condition: 
for each point of $X$ there is a neighbourhood
that is Sikorski diffeomorphic to a convex subset of a cartesian space.
Then the Sikorski structure on $X$ is reflexive. 
\end{theorem} 

\begin{proof}
By hypothesis, there is a collection $\{ V_k \}$ of subsets of $X$
whose interiors $\mathring{V}_k$ cover $X$
such that each $V_k$ is Sikorski diffeomorphic to a convex subset $C_k$
of a cartesian space.
Because a Sikorski diffeomorphism is in particular a homeomorphism,
each $\mathring{V}_k$ is diffeomorphic to an open subset $C_k'$ of $C_k$.
Because $X$ is locally compact, so is each $\mathring{V}_k$, 
and hence so is each $C'_k$.
By Lemma~\ref{l:locally convex}, each $C'_k$ is reflexive,
and hence so is each $\mathring{V}_k$.
By Proposition~\ref{p:local reflexivity}, $X$ is reflexive.
\end{proof}

\begin{remark}\labell{r:locally convex}
Theorem~\ref{t:locally convex} implies the following: let $K$ be a
locally closed convex set, and $F\colon K\to\RR^n$ a map.  Then $F$
is diffeologically smooth if and only if it is Sikorski smooth.  
\end{remark}

\begin{example}\labell{x:mfld-w-corners}
Let $M$ be a manifold-with-corners, and let $\calF$ and $\calD$ be its
natural Sikorski structure and diffeology, resp.  Then $(\calD,M,\calF)$
is a reflexive triple.  Indeed, recall that a manifold-with-corners of
dimension $n$ is by definition locally diffeomorphic to open subsets
of $[0,\infty)^n$.  
The result then follows from Theorem~\ref{t:locally convex}.
\end{example}



\section{Chen spaces}\labell{sec:chen}

Kuo-Tsai Chen introduced structures that are similar
to diffeologies but in which the domains of plots are convex
subsets of cartesian spaces rather than open subsets of cartesian spaces.
Chen's papers \cite{chen1973,chen1975, chen1977,chen1986}, 
which range over the years 1973 to 1986,
contain several different definitions of such spaces.
They converge to a definition \cite{chen1977,chen1986}
that is now typically taken to be ``the'' definition of a Chen space.
Baez and Hoffnung \cite[page 5793]{BH} identify these Chen spaces 
with concrete sheaves over the site of convex subsets of cartesian spaces.
Stacey's paper \cite[Section~3]{stacey}
compares different categories of smooth structures,
including Chen's different definitions (see \cite[p.~80]{stacey})
and diffeology, and describes functors between them.

Chen's smooth compatibility axiom
relies on the notion of a ``smooth map between convex sets'',
which  Chen does not explicitly define.
In practice, Chen seems to only need convex sets that are locally closed.
For these, there is an unambiguous notion 
of ``smooth map between convex sets'' (see Remark~\ref{r:locally convex} and Proposition~\ref{extend to boundary}).
However, Chen does not assume ``locally closed'' 
in his latest definitions \cite{chen1977,chen1986},
and without this assumption, 
as we have shown in \S\ref{s:example}, 
the notion of ``smooth map between convex sets'' is ambiguous.

We will interpret ``smooth map between convex sets''
as a map between the convex sets that is smooth 
\emph{in the diffeological sense}.
See Remark~\ref{r:equiv} for the rationale behind this choice.
For maps from a closed interval to $\R$,
this is equivalent to smoothness in the interior of the interval 
and one-sided smoothness at the boundary points.

Recall that a \emph{convex set} is a convex subset of a cartesian space,
and that the D-topology on a convex set coincides with its
subspace topology that is induced from its ambient cartesian space.

A \textbf{Chen structure} on a set $X$ is a collection $\calC_X$ of maps,
called \textbf{Chen plots}, from convex sets to $X$,
that satisfies the following axioms.
\begin{itemize}[itemsep=6pt]
\item[(C1)]
Constant maps to $X$ are Chen plots.
\item[(C2)]
The precomposition of a Chen plot
with a (diffeologically) smooth map between convex sets is a Chen plot.
\item[(C3)]
Given a convex set $C$ and a map $p \colon C \to X$,
if each point of $C$ has a convex neighbourhood $V$ in $C$ 
such that $p|_V$ is a Chen plot, then $p$ is a Chen plot.
\end{itemize}
A \textbf{Chen space} is a set equipped with a Chen structure.
A map from a Chen space $X$ to a Chen space $Y$ is \textbf{Chen smooth}
if its precomposition with every Chen plot of $X$ is a Chen plot of $Y$,
and it is a \textbf{Chen diffeomorphism} if it is bijective
and both it and its inverse are Chen smooth.
Chen spaces with Chen smooth maps form a category, which we denote $\Chen$.

Except where we say otherwise, 
we equip any subset of a cartesian space with the
\textbf{standard Chen structure}, consisting of the collection of 
all diffeologically smooth maps from convex sets to that subset. 

\begin{examples}\labell{x:chen}
A map between open subsets of cartesian spaces is Chen smooth 
if and only if it is smooth in the usual sense.
As a second example, 
the function $f\colon X\to\RR$ that we constructed in Section~\ref{s:example}
is Chen smooth.
\end{examples}

There is a natural functor $\Di\colon\Chen\to\Diffeol$ defined 
as follows:

\begin{itemize}[itemsep=6pt]
\item For any Chen space $(X,\calC_X)$, define $\Di(X,\calC_X)$ to be
the diffeological space $(X,\Di(\calC_X))$, 
in which a map $p\colon U\to X$ 
from an open subset of a Cartesian space
is in $\Di(\calC_X)$ if it is Chen smooth.
\item For any Chen smooth map $F\colon(X,\calC_X)\to(Y,\calC_Y)$, 
define $\Di(F)$ to be $F$.
\end{itemize}

Similarly, there is a natural function $\Ch\colon\Diffeol\to\Chen$ 
defined as follows:
\begin{itemize}[itemsep=8pt]
\item For any diffeological space $(X,\calD_X)$, define $\Ch(X,\calD_X)$
to be the Chen space $(X,\Ch(\calD_X))$, in which a map $p\colon C\to X$
from a convex set is in $\Ch(\calD_X)$ if it is diffeologically smooth.
\item For any diffeologically smooth map
$F\colon(X,\calD_X)\to(Y,\calD_Y)$, define $\Ch(F)$ to be $F$.
\end{itemize}

{
The functors $\Di$ and $\Ch$ are exactly the functors
$\mathfrak{So}$ and $\Ch^\sharp$ of Stacey's paper \cite{stacey}.
We refer the reader to \cite{stacey} 
for the proofs that $\Di$ indeed takes Chen structures
to diffeologies and Chen smooth maps to diffeologically smooth maps,
and that $\Ch$ indeed takes diffeologies to Chen structures
and diffeologically smooth maps to Chen smooth maps.
Moreover, Stacey shows that these functors are adjoint.  
A natural question is whether these functors
form an equivalence of categories between $\Chen$ and $\Diffeol$.  
They don't;\footnote{
In fact, according to Stacey, \cite[Corollary 8.9]{stacey},
there does not exist \emph{any} equivalence of categories
between $\Chen$ and $\Diffeol$.} see Example~\ref{nonstandard interval}.
The main issue is that there can be ``missing Chen plots''.
We will now show that, if we add these ``missing Chen plots'',
then we do obtain equivalent (in fact, isomorphic) categories.  

We define a Chen space $(X,\calC_X)$ to be \textbf{exhaustive} 
if it has the following property:
\begin{quotation}
Let $p \colon C \to X$ be a map from a convex set $C$ to $X$.
Assume that, for every convex open subset $U$ of a cartesian space
and smooth map $f \colon U \to C$,
the composition $p\circ f$ is in $\calC_X$. Then $p\in\calC_X$.
\end{quotation}
}

It is straightforward to show that the quality of being exhaustive is
invariant under Chen diffeomorphism. The following example is due to
Stacey \cite[page 100]{stacey}.

\begin{example} \labell{nonstandard interval}
Equip the interval $[0,1]$ with the set of those maps $p\colon C\to[0,1]$
from a convex set $C$ such that $p$ locally factors through a usual
smooth map $U\to[0,1]$, where $U$ is an open subset of a cartesian space.
This is a non-standard Chen structure on $[0,1]$.  It is not exhaustive;
for example, the identity map on $[0,1]$ is not in this Chen structure,
but its precomposition with any smooth map 
from an open convex subset $U$ of a cartesian space to $[0,1]$
is in this non-standard Chen structure. 
In contrast, the standard Chen structure on $[0,1]$ is exhaustive.
Thus, these two Chen structures on $[0,1]$ are not isomorphic.
But the functor $\Di$ takes both of these structures
to the standard diffeology on $[0,1]$.
So the functor $\Di$ is not one-to-one on isomorphism classes,
and so it is not an equivalence of categories.
\end{example}

Given a Chen space $(X,\calC_X)$, we define the \textbf{exhaustion} 
of $\calC_X$ to be the Chen structure $E(\calC_X):=\Ch(\Di(\calC_X))$ on~$X$.

\begin{lemma}[Properties of Exhaustions]\labell{l:exhaustion}
Let $(X,\calC_X)$ be a Chen space.
\begin{enumerate}
\item\labell{i:E makes exhaustive} $(X,E(\calC_X))$ is an exhaustive Chen space.
\item\labell{i:EC contains C} $\calC_X\subseteq E(\calC_X)$.
\item\labell{i:EEC is EC} If $(X,\calC_X)$ is exhaustive, 
then $E(\calC_X)=\calC_X$.
\end{enumerate}
\end{lemma}

The proofs rely on some results that we defer to an appendix.

\begin{lemma} \labell{l:Ch}
Any object in the image of the functor $\Ch$ is an exhaustive Chen space.
\end{lemma}

\begin{proof}
Fix a diffeological space $(X,\calD_X)$,
and consider $\calC_X := \Ch(\calD_X)$.
By Item \eqref{v} of Proposition \ref{p:DiCh}, $\Di(\Ch(\calD_X)) = \calD_X$.
Applying $\Ch$, we obtain $\Ch(\Di(\Ch(\calD_X))) = \Ch(\calD_X)$,
that is, $\Ch(\Di(\calC_X)) = \calC_X$.
By Item \eqref{ix} of Proposition \ref{p:DiCh},
this implies that $\calC_X$ is exhaustive.
\end{proof}

\begin{proof}[Proof of Lemma \ref{l:exhaustion}]
Item~\eqref{i:E makes exhaustive} is a consequence of Lemma~\ref{l:Ch}.

Item~\eqref{i:EC contains C} is Item \eqref{ix-2} of Proposition~\ref{p:DiCh},
but here is also a direct proof.
Let $p\colon C\to X$ be a Chen plot in $\calC_X$.  
Then $p$ is Chen smooth with respect to $\calC_X$ (Item~(\ref{r:Cplots})
of Lemma \ref{l:sanity}).
Because $\Di$ is a functor, $p$ is diffeologically smooth 
with respect to $\Di(\calC_X)$.
Because $\Ch$ is a functor, $p$ is Chen smooth 
with respect to $\Ch(\Di(\calC_X))$. 
So $p$ is a Chen plot in $\Ch(\Di(\calC_X))$ 
(Item~\eqref{r:Cplots} of Lemma \ref{l:sanity}).

To prove Item~\eqref{i:EEC is EC}, suppose $(X,\calC_X)$ is exhaustive. 
By Item~\eqref{i:EC contains C},
it suffices to show that $E(\calC_X)\subseteq \calC_X$.
This, in turns, follows from Item \eqref{ix} of Proposition~\ref{p:DiCh}.
\end{proof}

\begin{theorem}\labell{t:equiv}
The restriction of $\Di$ to exhaustive Chen spaces is an isomorphism 
of categories onto $\Diffeol$, with inverse $\Ch$.
\end{theorem}

\begin{proof}
By Lemma \ref{l:Ch}, $\Ch$ takes diffeological spaces
to exaustive Chen spaces.
By Item \eqref{v} of Proposition \ref{p:DiCh}, 
the composition $\Di \circ \Ch$
is the identity functor on the category of diffeological spaces.
By Part \eqref{i:EEC is EC} of Lemma \ref{l:exhaustion},
the composition $\Ch \circ \Di$ restricts to the identity functor
on the category of exhaustive Chen spaces.
\end{proof}

\begin{remark}\labell{r:equiv}
One motive for taking the notion of a ``smooth map between convex
subsets'' to be the diffeological notion is that an equivalent notion
is used in Baez and Hoffnung \cite{BH}.  A second motive is that this definition is equivalent to many other standard notions; see Proposition~\ref{extend to boundary} and \cite[page 5794]{BH}.
We now have a third motive:
the isomorphism of Theorem~\ref{t:equiv} no longer holds if we take
``smooth map between convex subsets'' to be the Sikorski notion.  Indeed,
the map $f\colon X\to\RR$ of Section~\ref{s:example} is Chen smooth when we take smoothness between convex subsets 
to be the diffeological notion.  However, if in the definition of ``Chen space''
we were to take smoothness between convex subsets to be the Sikorski notion,
then the map $f \colon X \to \R$ of Section~\ref{s:example}
would not be Chen smooth,
and hence the functor $\Di$ would no longer be full.
\end{remark}

\appendix
\section{Back and forth between Diffeologies and Chen structures}

We begin with a ``sanity check'' for our definitions.

\begin{lemma} \labell{l:sanity}
\begin{enumerate}[itemsep=6pt]
\item \labell{r:plots}
Given a diffeological space $(X,\calD_X)$ and a map $p \colon U \to X$
from an open subset $U$ of a cartesian space,
the map $p$ is diffeologically smooth if and only if it is a plot.
\item \labell{r:Cplots}
Given a Chen space $(X, \calC_X)$ and a map $p \colon C \to X$
from a convex set $C$, the map $p$
is Chen smooth if and only if it is a Chen plot in $\calC_X$.
\item \labell{r:sanity-S}
Given a Sikorski differential space $(X,\calF_X)$ and a function 
$f \colon X \to \R$, the function $f$ is Sikorski smooth
if and only if it is in $\calF_X$. 
\end{enumerate}
\end{lemma}

\begin{proof} \ 
We prove Item \eqref{r:plots}.
Every plot in $\calD_X$ is diffeologically smooth 
because $\calD_X$ satisfies the smooth compatibility axiom in the
definition of a diffeology. 
Every diffeologically smooth map 
from an open subset $U$ of a cartesian space to $X$ is a plot
by the definition of ``diffeologically smooth''
and because the identity map of $U$ is a plot.

We prove Item \eqref{r:Cplots}.
Every plot in $\calC_X$ is Chen smooth
because $\calC_X$ satisfies the smooth compatibility axiom
of a Chen structure. 
Every Chen smooth map from a convex subset $C$ to $X$
is a Chen plot by the definition of ``Chen smooth'' and because
the identity map on $C$ is a Chen plot.

We prove Item \eqref{r:sanity-S}.
Every function in $\calF_X$ is Sikorski smooth
because $\calF_X$ satisfies the smooth compatibility axiom in the definition
of a Sikorski structure. 
Every Sikorski smooth function
from $X$ to $\R$ is in $\calF_X$ by the definition of ``Sikorski smooth"
and because the identity map on $\R$ is in $\calF_{\R}$.
\end{proof}

We continue with a simple characterization of diffeological smoothness.

\begin{lemma} \labell{C to X}
Let $C$ and $X$ be diffeological spaces,
and let $p \colon C \to X$ be any map of their underlying sets.
Suppose that, for every plot $f \colon U \to C$ of $C$
whose domain $U$ is convex,
the composition $p \circ f \colon U \to X$ is a plot of $X$.
Then $p$ is diffeologically smooth.
\end{lemma} 

\begin{proof} 
Fix a plot $q \colon V \to C$ of $C$.
We need to prove that 
the composition $p \circ q \colon V \to X$ is a plot of $X$.
By the smooth compatibility axiom for the diffeology of $C$,
for every convex open subset $U$ of $V$, the restriction
$q|_U \colon U \to C$ is also a plot of $C$.
By the hypothesis on $p$, for every such $U$,
the composition $p \circ q|_U \colon U \to X$
is a plot of $X$. 
Because $V$ is covered by such convex open subsets $U$, 
and by the locality axiom for the diffeology of $X$,
\ $p\circ q \colon V \to X$ is a plot of $X$, as required.
\end{proof}

We now give the main result of this appendix.

\begin{proposition}\labell{p:DiCh}
\begin{enumerate}
\item \labell{v}
For any diffeological space $(X,\calD_X)$, 
$$\Di(\Ch(\calD_X)) = \calD_X.$$
\item \labell{ix}
For any Chen space $(X,\calC_X)$,
\ $\Ch(\Di(\calC_X))$ is the collection of those maps 
$p \colon C \to X$ from convex sets $C$ to $X$
that have the following property.
\begin{quotation}
For every convex open subset ${V}$ of a cartesian space
and every smooth map $f \colon {V} \to C$,
the composition $p \circ f \colon {V} \to X$ is a Chen plot in $\calC_X$.
\end{quotation}
Consequently, a Chen space $(X,\calC_X)$ is exhaustive
if and only if 
$$\Ch(\Di(\calC_X)) \subset \calC_X.$$
\item \labell{ix-2}
For any Chen space $(X,\calC_X)$,
$$\calC_X \subset \Ch(\Di(\calC_X)).$$
\end{enumerate}
\end{proposition}

\begin{proof} [Proof of Item \eqref{v} of Proposition~\ref{p:DiCh}]

Fix a diffeology $\calD_X$.

We will first prove that $\Di(\Ch(\calD_X)) \subset \calD_X$.
Take a plot $p \colon U \to X$ in $\Di(\Ch(\calD_X))$.
Because $\Di(\Ch(\calD_X))$ is a diffeology,
it satisfies the smooth compatibility axiom
in the definition of a diffeology,
so for any convex open subset ${V} \subset U$, the restriction 
$p|_{V} \colon {V} \to X$ is also a plot in $\Di(\Ch(\calD_X))$.
By the definition of $\Di(\Ch(\calD_X))$, 
every such $p|_{V} \colon {V} \to X$
is Chen smooth with respect to $\Ch(\calD_X)$.
By Item~(\ref{r:Cplots}) of Lemma~\ref{l:sanity}, every such $p|_{V}$ 
is a Chen plot in $\Ch(\calD_X)$.
By the definition of $\Ch(\calD_X)$, every such $p|_{V}$ 
is diffeologically smooth with respect to $\calD_X$.
By Item~(\ref{r:plots}) of Lemma~\ref{l:sanity}, Every such $p|_{{V}}$ is a plot in $\calD_X$.
Because $U$ can be covered by convex open subsets ${V}$,
and by the locality property of $\calD_X$, \ $p$ is a plot in $\calD_X$, 
as required.

We will now prove that $\calD_X \subset \Di(\Ch(\calD_X))$.
Take a plot $p \colon U \to X$ in $\calD_X$.
Because $\calD_X$ satisfies the smooth compatibility axiom
in the definition of a diffeology,
for every convex open subset ${V} \subset U$, 
the restriction $p|_{{V}} \colon {V} \to X$ is a plot in $\calD_X$.
By Item~(\ref{r:plots}) of Lemma~\ref{l:sanity}, every such $p|_{{V}}$ is diffeologically smooth
with respect to $\calD_X$.
Because $\Ch$ is a functor, every such $p|_{{V}}$ is Chen smooth
with respect to $\Ch(\calD_X)$.
Because $\Di$ is a functor, every such $p|_{{V}}$ is diffeologically smooth
with respect to $\Di(\Ch(\calD_X))$.
By Item~(\ref{r:plots}) of Lemma~\ref{l:sanity}, every such $p|_{{V}}$ 
is a plot in $\Di(\Ch(\calD_X))$.
Because $U$ can be covered by convex open subsets ${V}$,
and because $\Di(\Ch(\calD_X))$ 
satisfies the locality axiom in the definition of a diffeology,
$p$ is a plot in $\Di(\Ch(\calD_X))$, as required.
\end{proof}

\begin{proof} [Proof of Item \eqref{ix} of Proposition~\ref{p:DiCh}]
Fix a Chen structure $\calC_X$ and a map $p \colon C \to X$
from a convex set $C$ to $X$.

First, suppose that $p$ is a Chen plot in $\Ch(\Di(\calC_X))$,
and fix a convex open subset ${V}$ of a cartesian space
and a smooth map $f \colon {V} \to C$.
By the definition of $\Ch(\Di(\calC_X))$, \ $p$ is diffeologically smooth
with respect to $\Di(\calC_X)$.
Because the composition of diffeologically smooth maps 
is diffeologically smooth, the composition $p \circ f \colon {V} \to X$
is diffeologically smooth with respect to $\Di(\calC_X)$.
By Item~(\ref{r:plots}) of Lemma~\ref{l:sanity}, \ $p \circ f$ is a plot in $\Di(\calC_X)$.
By the definition of $\Di(\calC_X)$, 
\ $p \circ f$ is Chen smooth with respect to $\calC_X$.
By Item~(\ref{r:Cplots}) of Lemma~\ref{l:sanity}, \ $p \circ f$ is a Chen plot in $\calC_X$.

Next, suppose that for every convex open subset ${V}$ of a cartesian space
and every smooth map $f \colon {V} \to C$, the composition 
$p \circ f \colon {V} \to X$ is a Chen plot in $\calC_X$.
By Item~(\ref{r:Cplots}) of Lemma~\ref{l:sanity}, every such $p \circ f$ is Chen smooth
with respect to $\calC_X$.
By the definition of $\Di(\calC_X)$, every such $p \circ f$
is a plot in $\Di(\calC_X)$.
By Lemma~\ref{C to X},
$p$ is diffeologically smooth with respect to $\Di(\calC_X)$.
By the definition of $\Ch(\Di(\calC_X))$, 
$p$ is in $\Ch(\Di(\calC_X))$.
\end{proof}

\begin{proof} [Proof of Item \eqref{ix-2} of Proposition~\ref{p:DiCh}]
Fix a Chen space $(X,\calC_X)$.
Let $p \colon C \to X$ be a Chen plot in $\calC_X$.
Because $\calC_X$ satisfies the smooth compatibility axiom 
in the definition of a Chen structure, 
$p$ has the property listed in Item \eqref{ix}.
By Item~\eqref{ix}, \ $p$ is in $\Ch(\Di(\calC_X))$.
\end{proof}


\begin{thebibliography}{99}
\bibliographystyle{amsalpha}

\bibitem{BH}
John C.\ Baez and Alexander E.\ Hoffnung,
``Convenient categories of smooth spaces'',
\emph{Trans. Amer. Math. Soc.}, \textbf{363} (2011), no. 11, 5789--5825.

\bibitem{BIZKW}
Augustin T.\ Batubenge, Patrick Iglesias-Zemmour, Yael Karshon, 
and Jordan Watts, ''Diffeological, Fr\"{o}licher, and differential spaces''.\\
\url{https://arxiv.org/abs/1712.04576}

\bibitem{bierstone:diff fcts}
Edward Bierstone, ``Differentiable functions'', 
\emph{Bol.\ Soc.\ Brasil.\ Mat.} \textbf{11} (1980), no.\ 2, 139-189.

\bibitem{boman}
Jan Boman, ``Differentiability of a function and of its compositions
with functions of one variable'', \emph{Math. Scand.}, \textbf{20}
(1967), 249--268.

\bibitem{chen1973}
Kuo-Tsai Chen,
``Iterated integrals of differential forms and loop space homology'',
\emph{Ann.  of Math.} \textbf{97} (1973), no. 2, 217--246.

\bibitem{chen1975}
\rule{1cm}{0.01cm},
``Iterated integrals, fundamental groups and covering spaces'',
\emph{Trans.\ Amer.\ Math.\ Soc.} \textbf{206} (1975), 83--98.

\bibitem{chen1977}
\rule{1cm}{0.01cm}, ``Iterated path integral'',
\emph{Bull.\ of Amer.\ Math.\ Soc.} \textbf{83} (1977), no.~5, 831--879.

\bibitem{chen1986}
\rule{1cm}{0.01cm}, ``On differentiable spaces'',
In: \emph{Categories in continuum physics} (Buffalo, N.Y. 1982),
Lecture Notes in Math., \textbf{1174} (1986), Springer-Verlag, Berlin, 38--42.

\bibitem{CSW}
J.\ Daniel Christensen, Gordon Sinnamon, and Enxin Wu, ``The $D$-topology
for diffeological spaces'', \emph{Pacific J.\ Math.} \textbf{272} (2014), 
no.\ 1, 87-110.

\bibitem{frolicher}
Alfred Fr\"olicher, ``Smooth structures'', In: \emph{Category Theory (Gummersbach, 1981), Lecture Notes in Math.}, \textbf{962} (1982), Springer, New York, 69--81.

\bibitem{FK}
Alfred Fr\"olicher and Andreas Kriegl,
\emph{Linear Spaces and Differentiation Theory},
Wiley-Interscience, 1988.

\bibitem{iglesias}
Patrick Iglesias-Zemmour, \emph{Diffeology},
Math.\ Surveys and Monographs,
Amer.\ Math.\ Soc., 2012.

\bibitem{kriegl}
Andreas Kriegl, ``Remarks on germs in infinite dimensions'',
\emph{Acta Math.\ Univ.\ Comenian.\ (N.S.)} \textbf{66} (1997), 
no.\ 1, 117-134.

\bibitem{sniatycki}
Jedrzej \'Sniatycki,
\emph{Differential Geometry of Singular Spaces and Reduction of Symmetry},
New Mathematical Monographs (no.~23),
Cambridge University Press,
2013.

\bibitem{stacey}
Andrew Stacey, ``Comparative smootheology'',
\emph{Theory Appl.\ Categ.} \textbf{25} (2011), no.~4, 64--117.

\bibitem{virgin}
Bryce Virgin, ``Watts spaces and smooth maps'', \emph{PUMP: A Journal of Undergraduate Research}, \textbf{4} (2021), 63--85.

\bibitem{watts}
Jordan Watts, \emph{Diffeologies, Differential Spaces,
and Symplectic Geometry}, Ph.D.\ Thesis (2012), University of Toronto, Canada.
\end{thebibliography}
\end{document}